\documentclass{article}%
\usepackage{amssymb}
\usepackage{amsfonts}
\usepackage{ifpdf}
\usepackage{srcltx}
\usepackage{url}
\usepackage{hyperref}
\usepackage{cite}
\usepackage{color}
\usepackage{multirow}
\usepackage{graphicx}
\usepackage{caption}
\usepackage{epstopdf}
\usepackage{letltxmacro}
\usepackage[cmex10]{amsmath}
\usepackage{amssymb}
\usepackage{algorithmic}
\usepackage{array}
\usepackage{eqparbox}
\usepackage[tight,footnotesize]{subfigure}
\usepackage{tabularx}
\usepackage{hyphenat}
\usepackage{amsmath}
\usepackage[left=0.5in, right=0.5in, bottom=0.6in,top=0.6in]{geometry}%
\usepackage{graphicx}
\providecommand{\U}[1]{\protect\rule{.1in}{.1in}}

\begin{document}

\author{ M. El-Morshedy$^{1}$, M. S. Eliwa$^{1,}$\thanks{Corresponding author:
mseliwa@mans.edu.eg} \ and H. Nagy$^{1}$\\$^{1}$Mathematics Department, Faculty of Science, Mansoura \\University,\ Mansoura 35516, Egypt.}
\title{Exponentiated Discrete Lindley Distribution: Properties and Applications}
\date{}
\maketitle

\begin{abstract}
In this article, the exponentiated discrete Lindley distribution is presented
and studied. Some important distributional properties are discussed. Using the
maximum likelihood method, estimation of the model parameters is investigated.
Furthermore, simulation study is performed to observe the performance of the
estimates. Finally, the model with two real data sets is examined.


\end{abstract}

\section{Introduction}

Statistical (lifetime) distributions are commonly applied to describe and
predict real world phenomena. Several classical distributions have been
extensively used over the past decades for modeling data in several fields
such as engineering, medicine, finance, biological and actuarial science.
Lindley distribution (LiD) is one of the most important lifetime
distributions, it has some nice properties to be used in lifetime data
analysis, especially in applications modeling stress-strength model (see,
Lindley (1958)). This distribution can be shown as a mixture of exponential
and gamma distributions. The random variable (RV) $Z$ is said to have LiD with
one scale parameter $a>0$, if the cumulative distribution function (CDF) and
the probability density function (PDF) are given by
\begin{equation}%
{\textstyle\prod}
(z;a)=1-e^{-az}\left(  1+\frac{az}{a+1}\right)  ;\ \ z>0,\label{A}%
\end{equation}
and \
\begin{equation}
\pi(z;a)=\frac{a^{2}}{1+a}(z+1)e^{-az};\ \ z>0,\label{B}%
\end{equation}
respectively. Due to its wide applicability in many areas, several works aimed
at extending the LiD become very important. See, Ghitany et al. (2008a, 2008b,
2011, 2013), Zakerzadeh and Dolati (2009), Mahmoudi and Zakerzadeh (2010),
Jodr\'{a} (2010), Nadarajah et al. (2011), Bakouch et al. (2012), Merovci
(2013), Shanker and Mishra (2013), Shanker et al. (2013), Merovci and Elbatal
(2014), Merovci and Sharma (2014), Liyanage and Pararai (2014), Pararai et al.
(2015), Sharma et al. (2015), Nedjar and Zeghdoudi (2016), Zeghdoudi and
Nedjar (2016, 2017), \"{O}zel et al. (2016), Elbatal et al. (2016), Altun et
al. (2017), Mahmoud (2018), Jehhan et al (2018), among others.

Furthermore, some discrete versions of the LiD have been presented in the
statistical literature because in several cases, lifetimes need to be recorded
on a discrete scale rather than on a continuous analogue. So, discretizing
continuous distributions has received much attention in the statistical
literature. See, Sankaran (1970), Ghitany and Al-Mutairi (2009)\textbf{,
}Calder\'{\i}n-Ojeda and G\'{o}mez-D\'{e}niz (2013),\textbf{\ }Bakouch et al.
(2014), Tanka and Srivastava (2014), Munindra et al. (2015), Kus et al.
(2018), among others.

Also, several discrete distributions have been presented in the literature.
See, Roy (2003, 2004), Inusah and Kozubowski (2006), Krishna and Pundir
(2009), G\'{o}mez-D\'{e}niz (2010), Bebbington et al. (2012), Nekoukhou et al.
(2013), Vahid and Hamid (2015), Alamatsaz et al. (2016), Chandrakant et al.
(2017), among others. \ \

Although there are a number of discrete distributions in the literature, there
is still a lot of space left to develop new discretized distribution that is
suitable under different conditions. So, in this article, we introduce a
flexible discrete distribution called, the exponentiated discrete Lindley
distribution (EDLiD), because the discrete Lindley distribution (DLiD) does
not supply enough flexibility for analysis different types of lifetime data.

The article is organized as follows. In Section 2, we introduce the EDLiD.
Different statistical properties are studied in Section 3 . The estimation of
the model parameters by maximum likelihood is performed in Section 4. In
Section 5, simulation study is presented. Moreover, two applications to real
data illustrate the potentiality of the EDLiD. Finally, Section 6 provides
some conclusions.

\section{The EDLiD}

G\'{o}mez-D\'{e}niz and Calder\'{\i}n-Ojeda (2011) introduced the DLiD. The RV
$Y$ is said to have DLiD with a parameter $0<a<1$ if the CDF and the
probability mass function (PMF) are given by%
\begin{equation}
W(y;a)=\frac{1-a^{y+1}+\left[  \left(  2+y\right)  a^{y+1}-1\right]  \log
a}{1-\log a};\ \ y\in%
\mathbb{N}
_{0}=\left\{  0,1,2,3,...\right\}  ,
\end{equation}
and%
\begin{equation}
w(y,a)=\frac{a^{y}}{1-\log a}[a\log a+(1-a)(1-\log a^{y+1})];\ \ y\in%
\mathbb{N}
_{0},
\end{equation}
respectively. In the context of lifetime distributions with CDF $W(y)$, the
most widely used generalization technique is the exponentiated-W. Using this
method, for $b>0$, the CDF of the exponentiated-W class is given by%
\begin{equation}
F(y;a,b)=\left[  W(y;a)\right]  ^{b},
\end{equation}
(see, Lehmann (1952)). Therefore, the RV $X$ \ is said to have EDLiD with
shape parameter $b$ and scale parameter $a$ if the CDF and the reliability
function are given by%

\begin{equation}
\ \ F(x;a,b)=\frac{\Lambda(x+1;a,b)}{\left(  1-\log a\right)  ^{b}%
};\text{\ \ }x\in%
\mathbb{N}
_{0},
\end{equation}
and%
\begin{equation}
R(x;a,b)=\frac{\left(  1-\log a\right)  ^{b}-\Lambda(x+1;a,b)}{\left(  1-\log
a\right)  ^{b}};\text{\ \ }x\in%
\mathbb{N}
_{0},
\end{equation}
respectively, where
\begin{equation}
\Lambda(x;a,b)=\left(  1-a^{x}+\left[  \left(  1+x\right)  a^{x}-1\right]
\log a\right)  ^{b}.
\end{equation}
Further, the PMF of the EDLiD is given by%
\begin{equation}
f(x;a,b)=\frac{1}{\left(  1-\log a\right)  ^{b}}\left[  \Lambda
(x+1;a,b)-\Lambda(x;a,b)\right]  \ ;\ \ x\in%
\mathbb{N}
_{0},
\end{equation}
where $f(x;a,b)=F(x+1;a,b)-F(x;a,b)$. Figure 1 shows the plots of the PMF for
various values of the model parameters.%
\[%
{\parbox[b]{5.5486in}{\begin{center}
\includegraphics[
height=4.2298in,
width=5.5486in
]%
{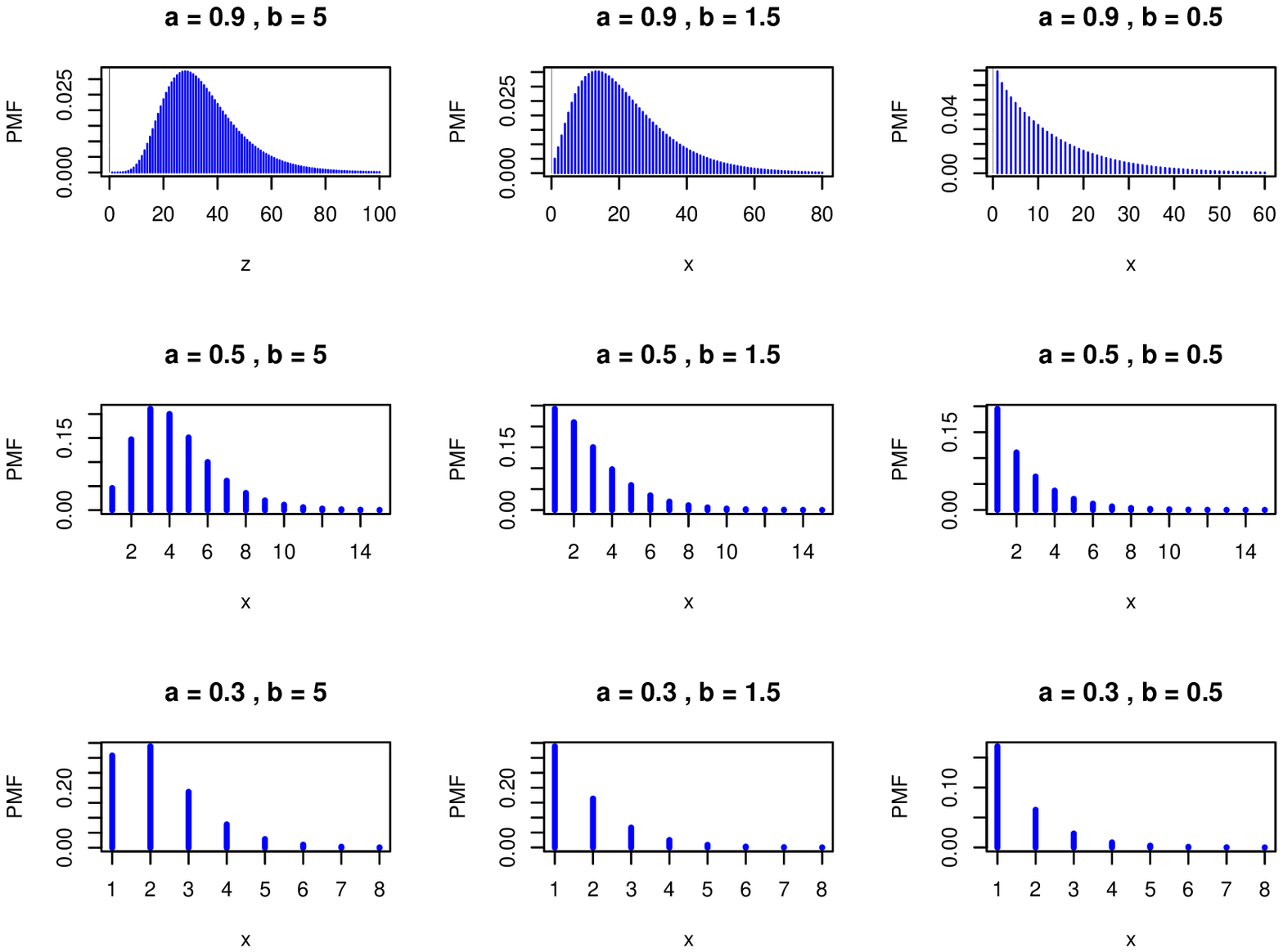}%
\\
Figure 1. The PMF of the EDLiD for various values of the parameters.
\end{center}}}
\]
From Figure 1, we note that the EDLiD can be take different shapes depending
on the values of the parameters. Moreover, the hazard rate function (Hrf) can
be expressed as%
\begin{equation}
h(x;a,b)=\frac{\Lambda(x+1;a,b)-\Lambda(x;a,b)}{\left(  1-\log a\right)
^{b}-\Lambda(x+1;a,b)};\text{\ \ }x\in%
\mathbb{N}
_{0}.
\end{equation}
Figure 2 shows the plots of the Hrf for various values of the model
parameters.%
\[%
{\parbox[b]{5.7493in}{\begin{center}
\includegraphics[
height=6.4524in,
width=5.7493in
]%
{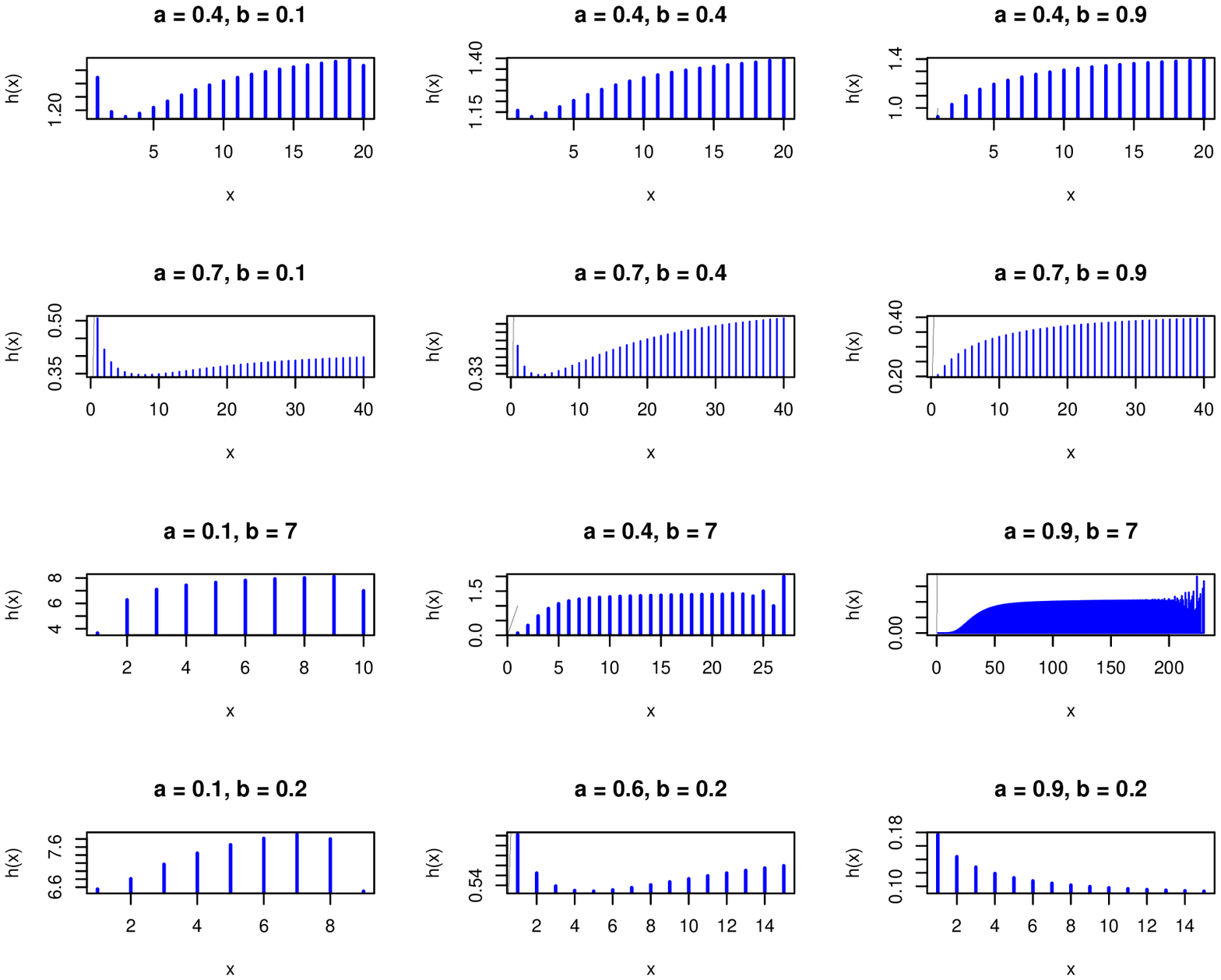}%
\\
Figure 2. The Hrf of the EDLiD for various values of the parameters.
\end{center}}}
\]

From Figure 2, it is clear that the Hrf can be increasing, decreasing, bathtub
and upside-down bathtub shaped. So, the EDLiD can be suitable for modeling
various data sets. Also, the reversed hazard rate function (Rhrf) of the EDLiD
can be expressed as follows%
\begin{equation}
r(x;a,b)=1-\frac{\Lambda(x;a,b)}{\Lambda(x+1;a,b)};\text{\ \ }x\in%
\mathbb{N}
_{0}.\text{ }%
\end{equation}
Figure 3 shows the plots of the Rhrf for various values of the model parameters.%

\[%
{\parbox[b]{5.047in}{\begin{center}
\includegraphics[
height=2.936in,
width=5.047in
]%
{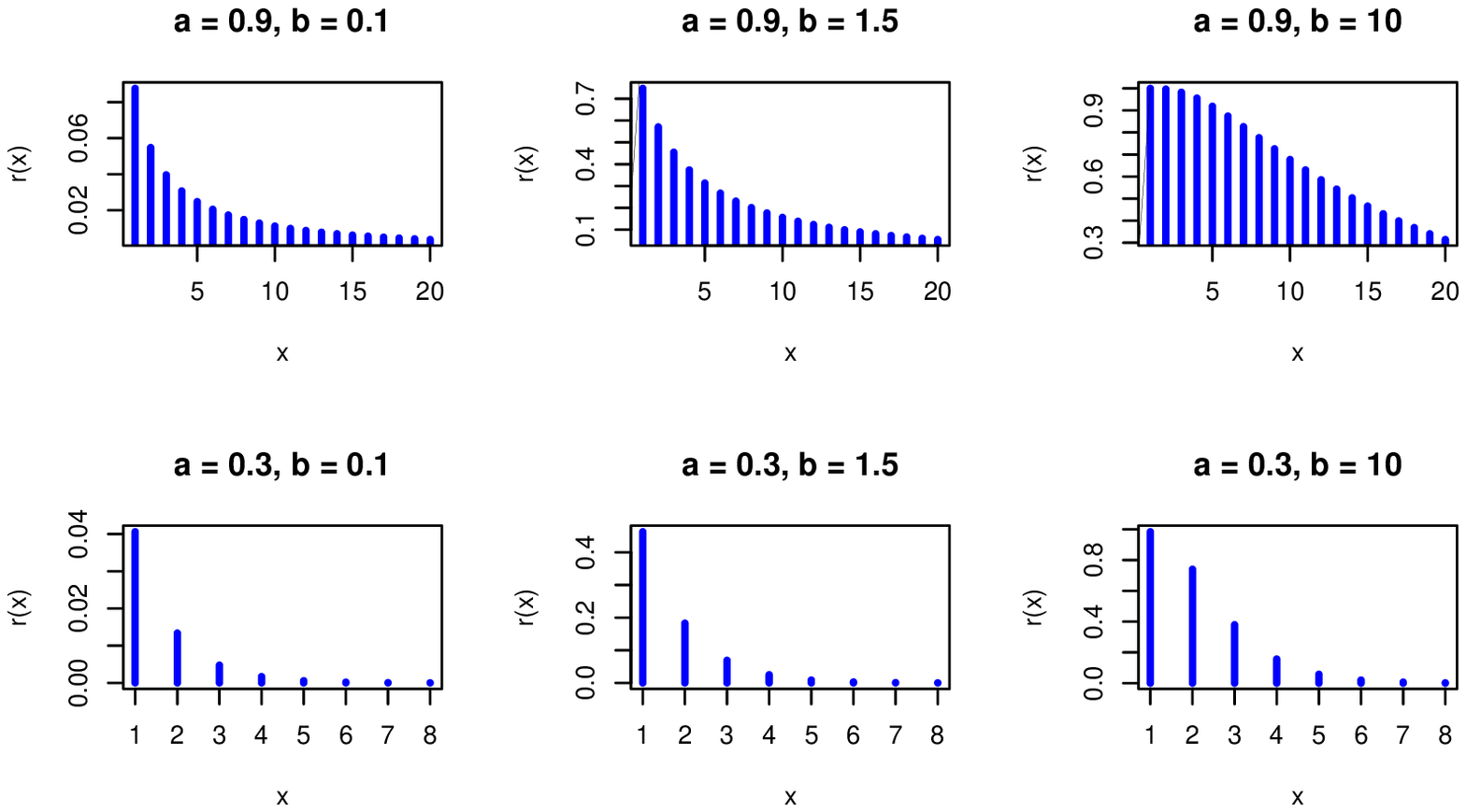}%
\\
Figure 3. The Rhrf of the EDLiD for various values of the parameters.
\end{center}}}
\]

\section{Different Properties}

\subsection{Moments}

Assume non-negative RV $X\sim$ EDLiD$(x;a,b)$. Then, the $r$th moment, say
$\varpi_{r}^{^{\prime}}$, is given by \ \ \
\begin{align}
\varpi_{r}^{^{\prime}}  & =\sum_{x=0}^{\infty}x^{r}f(x;a,b)=\sum_{x=1}%
^{\infty}\left[  x^{r}-\left(  x-1\right)  ^{r}\right]  R(x;a,b)\nonumber\\
& =\frac{1}{\left(  1-\log a\right)  ^{b}}\sum_{x=1}^{\infty}\left[
x^{r}-\left(  x-1\right)  ^{r}\right]  \left[  \left(  1-\log a\right)
^{b}-\Lambda(x+1;a,b)\right]  .\label{1-a}%
\end{align}
Using Equation (\ref{1-a}), we can get the mean ($\zeta$) and the variance
($\Upsilon$) of the random variable $X$ as follows
\begin{equation}
\zeta=\frac{1}{\left(  1-\log a\right)  ^{b}}\sum_{x=1}^{\infty}\left[
\left(  1-\log a\right)  ^{b}-\Lambda(x+1;a,b)\right]  ,
\end{equation}
and
\begin{equation}
\Upsilon=\frac{1}{\left(  1-\log a\right)  ^{b}}\sum_{x=1}^{\infty}\left[
2x-1\right]  \left[  \left(  1-\log a\right)  ^{b}-\Lambda(x+1;a,b)\right]
-\zeta^{2},
\end{equation}
respectively. Since $r$th moment is not in a closed form, then $\zeta$ and
$\Upsilon$ can only be numerically evaluated. Tables 1 and 2 obtain $\zeta$
and $\Upsilon$ of the EDLiD for different values of the model parameters respectively.%

\begin{align*}
& \text{Table 1.\ The }\zeta\ \text{of the EDLiD for different values of the
parameters }a\text{ and }b\text{.}\\
&
\begin{tabular}
[c]{|c|c|c|c|c|c|c|c|c|}\hline\hline
$b\downarrow a\rightarrow$ & $0.1$ & $0.2$ & $0.3$ & $0.4$ & $0.5$ & $0.6$ &
$0.7$ & $0.8$\\\hline\hline
$2$ & $0.364$ & $0.772$ & $1.269$ & $1.916$ & $2.816$ & $4.165$ & $6.424$ &
$10.969$\\\hline
$3$ & $0.508$ & $1.023$ & $1.626$ & $2.398$ & $3.463$ & $5.053$ & $7.708$ &
$13.033$\\\hline
$4$ & $0.631$ & $1.219$ & $1.893$ & $2.752$ & $3.934$ & $5.055$ & $8.637$ &
$14.526$\\\hline
$5$ & $0.737$ & $1.376$ & $2.102$ & $3.029$ & $4.305$ & $6.204$ & $9.365$ &
$15.693$\\\hline\hline
\end{tabular}
\end{align*}%
\begin{align*}
& \text{Table 2.\ The\ }\Upsilon\ \text{of the EDLiD for different values of
the parameters }a\text{ and }b\text{.}\\
&
\begin{tabular}
[c]{|c|c|c|c|c|c|c|c|c|}\hline\hline
$b\downarrow a\rightarrow$ & $0.1$ & $0.2$ & $0.3$ & $0.4$ & $0.5$ & $0.6$ &
$0.7$ & $0.8$\\\hline\hline
$2$ & $0.356$ & $0.793$ & $1.443$ & $2.514$ & $4.430$ & $8.225$ & $17.001$ &
$43.702$\\\hline
$3$ & $0.435$ & $0.875$ & $1.530$ & $2.623$ & $4.588$ & $8.479$ & $17.463$ &
$44.766$\\\hline
$4$ & $0.477$ & $0.899$ & $1.553$ & $2.624$ & $4.645$ & $8.568$ & $17.616$ &
$45.100$\\\hline
$5$ & $0.496$ & $0.901$ & $1.560$ & $2.675$ & $4.668$ & $8.599$ & $17.661$ &
$45.181$\\\hline\hline
\end{tabular}
\end{align*}
Depending on the model parameters, Tables 1 and 2 obtain that $\zeta$ and
$\Upsilon$ are increasing when $a$ is constant (increasing) and $b$ is
increasing (constant). Furthermore, the skewness ($\Xi$) and the kurtosis
($\Theta$) can be calculated as follows $\Xi=\frac{\varpi_{3}^{^{\prime}%
}-3\varpi_{2}^{^{\prime}}\varpi_{1}^{^{\prime}}+\varpi_{1}^{^{\prime}3}%
}{\Upsilon^{3/2}}$ and $\Theta=\frac{\varpi_{4}^{^{\prime}}-4\varpi
_{2}^{^{\prime}}\varpi_{1}^{^{\prime}}-3\varpi_{2}^{^{\prime}2}+12\varpi
_{2}^{^{\prime}}\mu\varpi_{1}^{^{\prime}2}-6\varpi_{1}^{^{\prime}4}}%
{\Upsilon^{2}}$. Tables 3 and 4 obtain the $\Xi$ and $\Theta$ of the EDLiD for
different values of the model parameters respectively.%
\begin{align*}
& \text{Table 3.\ The }\Xi\ \text{of the EDLiD for different values of the
parameters }a\text{ and }b\text{.}\\
&
\begin{tabular}
[c]{|c|c|c|c|c|c|c|c|c|}\hline\hline
$b\downarrow a\rightarrow$ & $0.1$ & $0.2$ & $0.3$ & $0.4$ & $0.5$ & $0.6$ &
$0.7$ & $0.8$\\\hline\hline
$2$ & \multicolumn{1}{|l|}{$1.667$} & \multicolumn{1}{|l|}{$1.335$} &
\multicolumn{1}{|l|}{$1.264$} & \multicolumn{1}{|l|}{$1.248$} &
\multicolumn{1}{|l|}{$1.2405$} & \multicolumn{1}{|l|}{$1.232$} &
\multicolumn{1}{|l|}{$1.223$} & \multicolumn{1}{|l|}{$1.215$}\\\hline
$3$ & \multicolumn{1}{|l|}{$1.222$} & \multicolumn{1}{|l|}{$1.076$} &
\multicolumn{1}{|l|}{$1.111$} & \multicolumn{1}{|l|}{$1.139$} &
\multicolumn{1}{|l|}{$1.148$} & \multicolumn{1}{|l|}{$1.148$} &
\multicolumn{1}{|l|}{$1.145$} & \multicolumn{1}{|l|}{$1.142$}\\\hline
$4$ & \multicolumn{1}{|l|}{$0.966$} & \multicolumn{1}{|l|}{$0.980$} &
\multicolumn{1}{|l|}{$1.066$} & \multicolumn{1}{|l|}{$1.098$} &
\multicolumn{1}{|l|}{$1.108$} & \multicolumn{1}{|l|}{$1.110$} &
\multicolumn{1}{|l|}{$1.109$} & \multicolumn{1}{|l|}{$1.107$}\\\hline
$5$ & \multicolumn{1}{|l|}{$0.809$} & \multicolumn{1}{|l|}{$0.957$} &
\multicolumn{1}{|l|}{$1.051$} & \multicolumn{1}{|l|}{$1.075$} &
\multicolumn{1}{|l|}{$1.084$} & \multicolumn{1}{|l|}{$1.088$} &
\multicolumn{1}{|l|}{$1.088$} & \multicolumn{1}{|l|}{$1.088$}\\\hline\hline
\end{tabular}
\end{align*}

\begin{align*}
& \text{Table 4.\ The }\Theta\ \text{of the EDLiD for different values of the
parameters }a\text{ and }b\text{.}\\
&
\begin{tabular}
[c]{|c|c|c|c|c|c|c|c|c|}\hline\hline
$b\downarrow a\rightarrow$ & $0.1$ & $0.2$ & $0.3$ & $0.4$ & $0.5$ & $0.6$ &
$0.7$ & $0.8$\\\hline\hline
$2$ & \multicolumn{1}{|l|}{$6.100$} & \multicolumn{1}{|l|}{$5.491$} &
\multicolumn{1}{|l|}{$5.466$} & \multicolumn{1}{|l|}{$5.468$} &
\multicolumn{1}{|l|}{$5.452$} & \multicolumn{1}{|l|}{$5.423$} &
\multicolumn{1}{|l|}{$5.393$} & \multicolumn{1}{|l|}{$5.368$}\\\hline
$3$ & \multicolumn{1}{|l|}{$4.679$} & \multicolumn{1}{|l|}{$4.923$} &
\multicolumn{1}{|l|}{$5.135$} & \multicolumn{1}{|l|}{$5.200$} &
\multicolumn{1}{|l|}{$5.210$} & \multicolumn{1}{|l|}{$5.201$} &
\multicolumn{1}{|l|}{$5.189$} & \multicolumn{1}{|l|}{$5.178$}\\\hline
$4$ & \multicolumn{1}{|l|}{$4.181$} & \multicolumn{1}{|l|}{$4.860$} &
\multicolumn{1}{|l|}{$5.051$} & \multicolumn{1}{|l|}{$5.097$} &
\multicolumn{1}{|l|}{$5.108$} & \multicolumn{1}{|l|}{$5.108$} &
\multicolumn{1}{|l|}{$5.104$} & \multicolumn{1}{|l|}{$5.099$}\\\hline
$5$ & \multicolumn{1}{|l|}{$4.051$} & \multicolumn{1}{|l|}{$4.845$} &
\multicolumn{1}{|l|}{$5.002$} & \multicolumn{1}{|l|}{$5.038$} &
\multicolumn{1}{|l|}{$5.055$} & \multicolumn{1}{|l|}{$5.061$} &
\multicolumn{1}{|l|}{$5.059$} & \multicolumn{1}{|l|}{$5.059$}\\\hline\hline
\end{tabular}
\end{align*}
Tables 3 and 4 obtain that $\Xi$ and $\Theta$ are decreasing when $a$ is
constant and $b$ is increasing. On the other hand, we can get the probability
generating function (PGF) of the RV $X$ as a form%
\begin{align}
\Omega_{X}(t)  & =\sum_{x=0}^{\infty}t^{x}f(x;a,b)\nonumber\\
& =R(0)+(t-1)R(1)+(t^{2}-t)R(2)+(t^{3}-t^{2})R(3)+...\nonumber\\
& =1+\frac{(t-1)}{\left(  1-\log a\right)  ^{b}}\sum_{x=1}^{\infty}%
t^{x-1}\left[  \left(  1-\log a\right)  ^{b}-\Lambda(x+1;a,b)\right]
.\label{JJJ}%
\end{align}
Using Equation (\ref{JJJ}), we can get $\zeta$ and $\Upsilon$ of the RV $X$ as
a form $\zeta=\frac{d}{dt}\Omega_{X}(t)|_{t=1}$ and $\Upsilon=$ $\frac{d^{2}%
}{dt^{2}}\Omega_{X}(t)|_{t=1}+\frac{d}{dt}\Omega_{X}(t)|_{t=1}-\left(
\frac{d}{dt}\Omega_{X}(t)|_{t=1}\right)  ^{2}$.

\subsection{Mean residual lifetime ($\varsigma(i)$) and mean past lifetime
($\varsigma^{\ast}(i)$)}

In order to study the ageing behavior of a component or a system of components
there have been defined several measures in the reliability and survival
analysis literature. The $\varsigma(i)$ is a helpful tool to model and analyze
the burn-in and maintenance policies. In the discrete setting, $\varsigma(i)$
is defined as%
\begin{equation}
\varsigma(i)=E\left(  T-i|T\geq i\right)  =\frac{1}{R(i)}\sum_{j=i+1}%
^{l}R(j)\ ;\ \ i\in%
\mathbb{N}
_{0},
\end{equation}
where $0<l<\infty$. If the RV $T\sim$ EDLiD$(a,b)$, then the $\varsigma(i) $
can be expressed as follows%
\begin{equation}
\varsigma(i)=\frac{1}{\left(  1-\log a\right)  ^{b}-\Lambda(i+1;a,b)}%
\sum_{j=i+1}^{l}\left[  \left(  1-\log a\right)  ^{b}-\Lambda(j+1;a,b)\right]
.
\end{equation}
Another measure of interest in survival analysis is $\varsigma^{\ast}(i)$. It
measures the time elapsed since the failure of $T$ given that the system has
failed sometime before $i$. In the discrete setting, $\varsigma^{\ast}(i)$ is
defined as%
\begin{equation}
\varsigma^{\ast}(i)=E\left(  i-T|T<i\right)  =\frac{1}{F(i-1)}\sum_{m=1}%
^{i}F(m-1)\ ;\ \ i\in%
\mathbb{N}
_{0}-\{0\},
\end{equation}
where $\varsigma^{\ast}(i)=0$ (see, Goliforushani and Asadi (2008)). If the RV
$T\sim$ EDLiD$(a,b)$, then the $\varsigma^{\ast}(i)$ can be represented as
follows%
\begin{equation}
\varsigma^{\ast}(i)=\frac{1}{\Lambda(i;a,b)}\sum_{m=1}^{i}\Lambda
(m;a,b).\label{GG}%
\end{equation}
For $i\in%
\mathbb{N}
_{0},$ we get $\varsigma^{\ast}(i)\leq i$.

\textbf{Lemma 1}. The mean of the RV $T\sim$ EDLiD$(a,b)$ can be expressed as%
\[
\zeta=i-\frac{\Lambda(i;a,b)}{\left(  1-\log a\right)  ^{b}}\varsigma^{\ast
}(i)+\frac{\left(  1-\log a\right)  ^{b}-\Lambda(i+1;a,b)}{\left(  1-\log
a\right)  ^{b}}\varsigma(i)\ ;\ \ i\in%
\mathbb{N}
_{0}.
\]

\textbf{Proof. }It is easy to prove this Lemma by using the following Equation%

\[
\zeta=\varsigma(0)=\sum_{j=1}^{l}R(j;a,b)=\sum_{j=1}^{i}R(j;a,b)+\sum
_{j=i+1}^{l}R(j;a,b).
\]

\textbf{Lemma 2. }The Rhrf and the $\varsigma^{\ast}(i;a,b)$ are related as
follows%
\begin{equation}
r(i;a,b)=\frac{1-\varsigma^{\ast}(i+1;a,b)+\varsigma^{\ast}(i;a,b)}%
{\varsigma^{\ast}(i;a,b)}\ ;\ \ i\in%
\mathbb{N}
_{0}-\{0\}.\
\end{equation}
\ \ \ \ \ \

\textbf{Proof. }%
\begin{align*}
F(i;a,b)\varsigma^{\ast}(i+1;a,b)-F(i-1;a,b)\varsigma^{\ast}(i;a,b)  &
=\sum_{j=1}^{i+1}F(j-1;a,b)-\sum_{j=1}^{i}F(j-1;a,b)\\
& =F(i;a,b).
\end{align*}
Dividing both sides of this Equation by $F(i;a,b)>0$, and noting that
$1-r(i;a,b)=\frac{F(i-1;a,b)}{F(i;a,b)}$, we get the required result.

\subsection{Stress-strength (S-S$^{\ast}$) analysis}

S-S$^{\ast}$ analysis has been used in mechanical component design. The
probability of failure is based on the probability of S exceeding S$^{\ast}$.
Assume that both S and S$^{\ast}$ are in the positive domain. The expected
reliability ($R^{\ast}$) can be calculated by
\begin{equation}
R^{\ast}=P\left[  X_{S}\leq X_{S^{\ast}}\right]  =\sum_{x=0}^{\infty}f_{X_{S}%
}(x)R_{X_{S^{\ast}}}(x)\text{ \ .}%
\end{equation}
If $X_{S}\sim$ EDLiD$(a_{1},b_{1})$ and $X_{S^{\ast}}\sim$ EDLiD$(a_{2}%
,b_{2})$, then
\begin{equation}
R^{\ast}=\frac{\sum_{x=0}^{\infty}\left[  \Lambda(x+1;a_{1},b_{1}%
)-\Lambda(x;a_{1},b_{1})\right]  \left[  \left(  1-\log a_{2}\right)  ^{b_{2}%
}-\Lambda(x+1;a_{2},b_{2})\right]  }{\left(  1-\log a_{1}\right)  ^{b_{1}%
}\left(  1-\log a_{2}\right)  ^{b_{2}}}.\label{100}%
\end{equation}
From Equation (\ref{100}), it is clear that the value of $R^{\ast}$ does not
depend only on the values of the model parameters.

\subsection{Order statistics (Os) and L-moment (Lm) statistics}

Let $X_{1},X_{2},...$,$X_{n}$ be a random sample from the EDLiD, and let
$X_{1:n},X_{2:n},...,X_{n:n}$ be their corresponding Os. Then, the CDF of
$i$th Os for an integer value of $x$ can be expressed as%
\begin{align}
F_{i:n}(x;a,b)  & =\sum_{k=i}^{n}\left(
\begin{array}
[c]{c}%
n\\
k
\end{array}
\right)  \left[  F_{i}(x;a,b)\right]  ^{k}\left[  1-F_{i}(x;a,b)\right]
^{n-k}\nonumber\\
& =\sum_{k=i}^{n}\sum_{j=0}^{n-k}\circleddash_{(j)}^{(n,k)}\frac
{\Lambda(x+1;a,b(k+j))}{\left(  1-\log a\right)  ^{b(k+j)}},
\end{align}
where $\circleddash_{(j)}^{(n,k)}=(-1)^{j}\left(
\begin{array}
[c]{c}%
n\\
k
\end{array}
\right)  \left(
\begin{array}
[c]{c}%
n-k\\
j
\end{array}
\right)  $. Furthermore, the PMF of the $i$th Os can be expressed as%
\begin{equation}
f_{i:n}(x;a,b)=\sum_{k=i}^{n}\sum_{j=0}^{n-k}\circleddash_{(j)}^{(n,k)}%
\frac{\left[  \Lambda(x+1;a,b(k+j))-\Lambda(x;a,b(k+j))\right]  }{\left(
1-\log a\right)  ^{b(k+j)}}.
\end{equation}
So, the $vth$ moments of $X_{i:n}$ can be written as%
\begin{equation}
E(X_{i:n}^{v})=\sum_{x=0}^{\infty}\sum_{k=i}^{n}\sum_{j=0}^{n-k}%
\circleddash_{(j)}^{(n,k)}x^{v}\frac{\left[  \Lambda(x+1;a,b(k+j))-\Lambda
(x;a,b(k+j))\right]  }{\left(  1-\log a\right)  ^{b(k+j)}}.
\end{equation}

On the other hand, Hosking (1990) has defined the L-moments (Lms) to summaries
theoretical distribution and observed samples. He has shown that the Lms have
good properties as measure of distributional shape and are useful for fitting
distribution to data. Lms are expectation of certain linear combinations of
Os. The Lms of \ $X$ \ can be expressed as follows
\begin{equation}
\Delta_{s}=\frac{1}{s}\sum_{j=0}^{s-1}(-1)^{j}\left(
\begin{array}
[c]{c}%
s-1\\
j
\end{array}
\right)  E\left(  X_{s-j:s}\right)  .\label{OO}%
\end{equation}
Since Hosking has defined the Lms of $X$ to be the quantities. Then, we can
propose some statistical measures such as Lm of mean $(Lm-M)=\Delta_{1}$, Lm
coefficient of variation $(Lm-Cv)=\frac{\Delta_{2}}{\Delta_{1}}$, Lm
coefficient of skewness $(Lm-Sk)=\frac{\Delta_{3}}{\Delta_{2}}$, and Lm
coefficient of kurtosis $(Lm-Ku)=\frac{\Delta_{4}}{\Delta_{2}}$.

\section{Estimation}

In this section, we determine the maximum likelihood estimates (MLEs) of the
model parameters from complete samples. Assume $X_{1},X_{2},...,X_{n}$ be a
random sample of size $n$ from the EDLiD($a,b$). The log-likelihood function
($L$) can be expressed as%
\begin{equation}
L(x;a,b)=-nb\log\left(  1-\log a\right)  +\sum_{i=1}^{n}\log\left[
\Lambda(x+1;a,b)-\Lambda(x;a,b)\right]  .\label{AA}%
\end{equation}
By differentiating Equation (\ref{AA}) with respect to the parameters $a$ and
$b$, we get the normal nonlinear likelihood equations as follows%
\begin{equation}
\frac{n\widehat{b}}{\widehat{a}\left(  1-\log\widehat{a}\right)  }+\widehat
{b}\sum_{i=1}^{n}\frac{\left[  V_{1}(x_{i}+1;\widehat{a})\right]
^{\widehat{b}-1}V_{2}(x_{i}+1;\widehat{a})-\left[  V_{1}(x_{i};\widehat
{a})\right]  ^{\widehat{b}-1}V_{2}(x_{i};\widehat{a})}{\Lambda(x_{i}%
+1;\widehat{a},\widehat{b})-\Lambda(x_{i};\widehat{a},\widehat{b})}=0,
\end{equation}
and%
\begin{equation}
-n\log\left(  1-\log\widehat{a}\right)  +\sum_{i=1}^{n}\frac{\Lambda
(x_{i}+1;\widehat{a},\widehat{b})\log(V_{1}(x_{i}+1;\widehat{a}))-\Lambda
(x_{i};\widehat{a},\widehat{b})\log(V_{1}(x_{i};\widehat{a}))}{\Lambda
(x_{i}+1;\widehat{a},\widehat{b})-\Lambda(x_{i};\widehat{a},\widehat{b})}=0,
\end{equation}
respectively, where $V_{1}(x;\widehat{a})=1-\widehat{a}^{x}+\left[  \left(
1+x\right)  \widehat{a}^{x}-1\right]  \log\widehat{a}$, and $V_{2}%
(x;\widehat{a})=x(x+1)\widehat{a}^{x-1}\log\widehat{a}-x\widehat{a}%
^{x-1}+\frac{1}{\widehat{a}}\left[  \left(  1+x\right)  \widehat{a}%
^{x}-1\right]  .$ Analytical software is required to get the values of the
model parameters.

\section{Applications}

\subsection{Simulation results}

In this section, we obtain the behavior of the MLEs of the EDLiD for a sample
size $n$ using a simulation study. At first, to generate a RV $X$ from the
EDLiD, we generate the value $z$ from the continuous ELiD. Then, discretize
this value to obtain $x$. The steps for a simulation study: choose the initial
values of the model parameters, say EDLiD($0.8,0.9$), generate $m=1000$
samples of size $n;n=$ 25, 50, 100, 150, 200, 250, 300, 350, 400, 450, 500,
550, 600, 650, 700, 750, compute the MLE's for the $m$ samples, say
{$(\widehat{a}_{i},\widehat{b}_{i}$}$)$ for $i=1,2,\dots,m$. Finally, compute
the average of biases and the average of mean squared errors (MSE(.)). Figure
4 shows how the biases and MSE vary with respect to $n$.%

\[%
{\parbox[b]{5.7493in}{\begin{center}
\includegraphics[
height=2.738in,
width=5.7493in
]%
{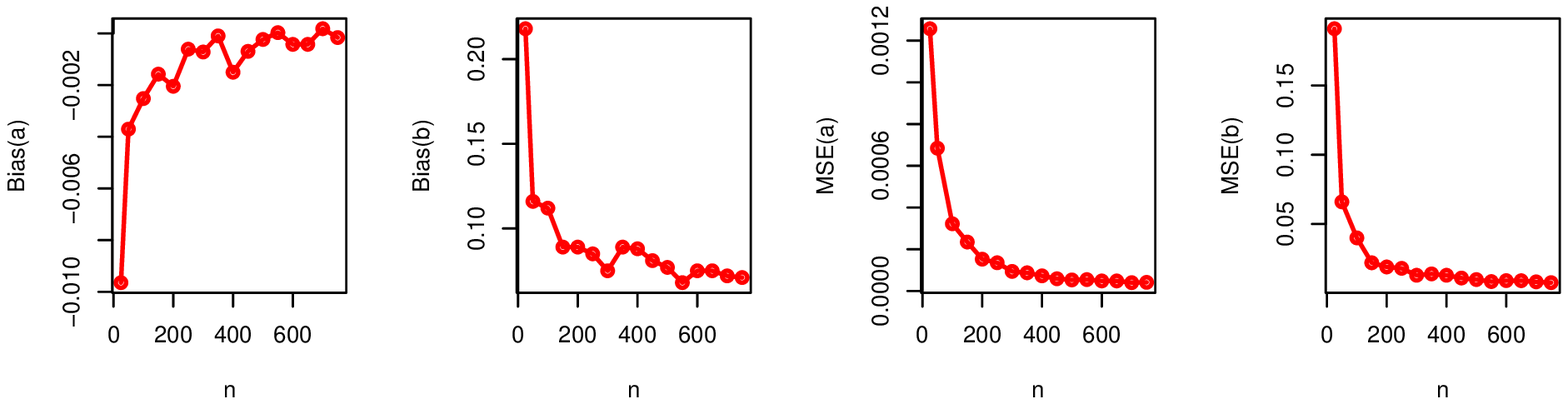}%
\\
Figure 4. The plots of bias(a), bias(b), MSE(a) and MSE(b) versus $n$ = 25,
50, 100,150, 200, ...,750.
\end{center}}}
\]
From Figure 4, it is clear that the biases and the MSEs of the estimated
parameters $\longrightarrow0$ while $n$ growing. So, the MLE is a good method
for estimating the model parameters.

\subsection{Data analysis}

In this section, we illustrate the importance of the EDLiD using two real data sets.

\textbf{The first data set (I): }represents the number of women who are
working on shells for 5 weeks discussed in Consul and Jain (1973). We shall
compare the fits of the EDLiD with some competitive models such as discrete
generalized exponential second type (DGE$_{2}$), discrete Weibull (DW),
discrete Lindley (DLi), discrete Pareto (DPa) and Poisson (P)
distributions\textbf{.}

\textbf{The second data set (II): }represents the counts of cysts of kidneys
using steroids. This data set originated from a study Chan et al. (2009). We
shall compare the fits of the EDLiD with some competitive models such as DW,
discrete Burr-XII (DB-XII), discrete Lomax (DLo), geometric (Geo), DLi, P and
discrete Rayleigh (DR) distributions.

The fitted models are compared using some criteria namely, the maximized
log-likelihood ($-L$), Akaike Information Criterion ($AIC$), Correct Akaike
Information Criterion ($CAIC$), Bayesian Information Criterion ($BIC$),
Hannan-Quinn Information Criterion ($HQIC$), chi-square ($\chi^{2}$) and its
P-value. \

\textbf{For the data set (I)},\textbf{\ }Tables 5 and 6 obtain the MLEs with
their corresponding standard errors (Se(.)), as well as $-L$, AIC, CAIC, BIC,
HQIC, $\chi^{2}$, degree of freedom (d.f), observed frequency (OF), expected
frequency (EF) and P-values respectively.%
\begin{align*}
& \text{Table 5. The MLEs with their corresponding Se for data set I.}\\
&
\begin{tabular}
[c]{|c|c|c|c|c|}\hline\hline
\textbf{Model }$\downarrow$ \textbf{Parameter} $\rightarrow$ & $\widehat{a} $
& \textbf{Se (}$\widehat{a}$\textbf{)} & $\widehat{b}$ & \textbf{Se
(}$\widehat{b}$\textbf{)}\\\hline\hline
\textbf{EDLi} & $0.263$ & $0.026$ & $0.693$ & $0.097$\\\hline
\textbf{DGE}$_{2}$ & $0.338$ & $0.031$ & $0.899$ & $0.125$\\\hline
\textbf{DW} & $0.312$ & $0.018$ & $0.967$ & $0.054$\\\hline
\textbf{DLi} & $0.209$ & $0.011$ & $-$ & $-$\\\hline
\textbf{DPa} & $0.139$ & $0.011$ & $-$ & $-$\\\hline
\textbf{P} & $0.466$ & $0.027$ & $-$ & $-$\\\hline\hline
\end{tabular}
\end{align*}

\begin{align*}
& \text{Table 6. The goodness of fit tests for data set I.}\\
&
\begin{tabular}
[c]{|c|c|c|c|c|c|c|c|}\hline\hline
$\mathbf{X}$ & \textbf{OF} & \multicolumn{6}{|c|}{\textbf{EF}}\\\cline{3-8}%
\cline{3-8}
&  & \textbf{EDLi} & \textbf{DGE}$_{2}$ & \textbf{DW} & \textbf{DLi} &
\textbf{DPa} & \textbf{P}\\\hline\hline
\textbf{0} & $447$ & $446.91$ & $446.71$ & $445.54$ & $430.06$ & $482.83$ &
$406.31$\\
\textbf{1} & $132$ & $131.87$ & $133.53$ & $135.36$ & $154.66$ & $90.57$ &
$189.03$\\
\textbf{2} & $42$ & $45.84$ & $44.29$ & $43.99$ & $45.75$ & $31.94$ &
$43.97$\\
\textbf{3} & $21$ & $15.28$ & $14.89$ & $14.62$ & $12.35$ & $14.87$ & $6.82$\\
\textbf{4} & $3$ & $4.91$ & $5.02$ & $4.92$ & $3.16$ & $8.11$ & $0.79$\\
$\mathbf{\geq5}$ & $2$ & $2.19$ & $2.56$ & $2.75$ & $1.02$ & $18.68$ & $0.08
$\\\hline
\textbf{Total} & $\mathbf{647}$ & $\mathbf{647}$ & $\mathbf{647}$ &
$\mathbf{647}$ & $\mathbf{647}$ & $\mathbf{647}$ & $\mathbf{647}$\\\hline
$\mathbf{-L}$ &  & $591.9$ & $592.2$ & $592.3$ & $595.3$ & $733.5$ & $617.2$\\
\textbf{AIC} &  & $1187.8$ & $1188.4$ & $1188.6$ & $1192.6$ & $1239.7$ &
$1236.4$\\
\textbf{CAIC} &  & $1187.8$ & $1188.4$ & $1188.6$ & $1192.6$ & $1239.7$ &
$1236.4$\\
\textbf{BIC} &  & $1196.8$ & $1197.3$ & $1197.5$ & $1197.0$ & $1244.2$ &
$1240.8$\\
\textbf{HQIC} &  & $1191.3$ & $1191.8$ & $1192.1$ & $1194.3$ & $1241.4$ &
$1238.1$\\\hline
$\mathbf{\chi}^{2}$ &  & $3.084$ & $3.521$ & $3.79$ & $10.514$ & $45.029$ &
$70.457$\\
\textbf{d.f} &  & $2$ & $2$ & $2$ & $3$ & $3$ & $3$\\
\textbf{P.value} &  & $0.21395$ & $0.1719$ & $0.1503$ & $0.01466$ & $<0.01$ &
$<0.01$\\\hline\hline
\end{tabular}
\end{align*}

From Table 6, it is clear that the EDLiD is the best distribution among all
tested distributions, because it has the smallest value among $-L$, AIC, CAIC,
BIC, HQIC and $\chi^{2}$, as well as it has the largest P-value. Figure 5
shows the fitted PMFs for data set I, which support the results in Table 6.%
\[%
{\parbox[b]{5.047in}{\begin{center}
\includegraphics[
height=4.0421in,
width=5.047in
]%
{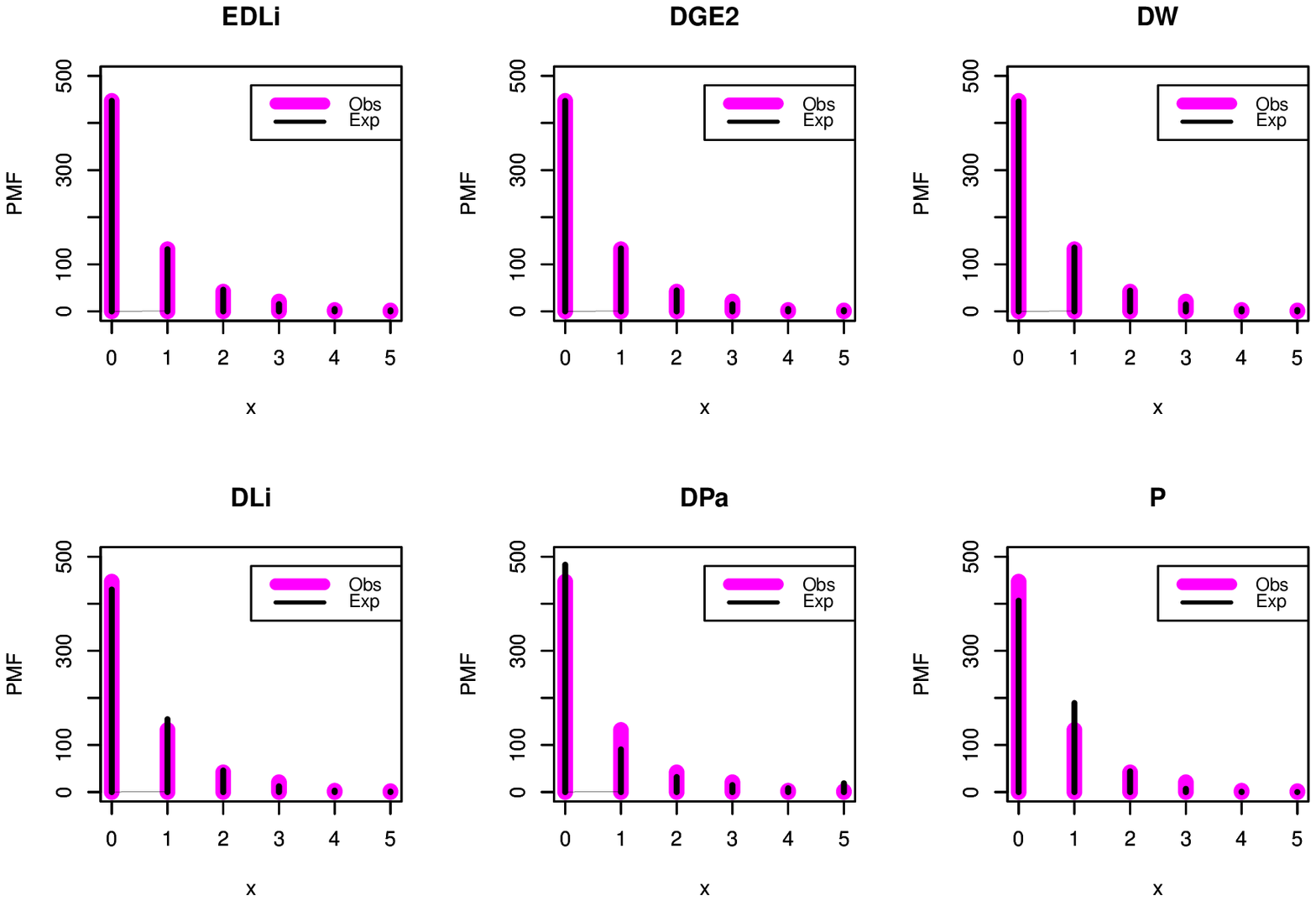}%
\\
Figure 5. The fitted PMFs for data set I.
\end{center}}}%
\]
\textbf{For the data set (II)},\textbf{\ }Tables 7 and 8 obtain the MLEs with
their corresponding Se, as well as $-L$, AIC, CAIC, BIC, HQIC, OF, EF, d.f,
$\chi^{2}$and P-values respectively.%

\begin{align*}
& \text{Table 7. The MLEs with their corresponding Se for data set II.}\\
&
\begin{tabular}
[c]{|c|c|c|c|c|c|c|}\hline\hline
\textbf{Model }$\downarrow$ \textbf{Parameter} $\rightarrow$ & $\widehat{a} $
& \textbf{Se (}$\widehat{a}$\textbf{)} & $\widehat{b}$ & \textbf{Se
(}$\widehat{b}$\textbf{)} & $\widehat{c}$ & \textbf{Se (}$\widehat{c}%
$\textbf{)}\\\hline\hline
\textbf{EDLi} & $0.672$ & $0.048$ & $0.264$ & $0.056$ & $-$ & $-$\\\hline
\textbf{DW} & $0.421$ & $0.047$ & $0.629$ & $0.073$ & $-$ & $-$\\\hline
\textbf{DB-XII} & $0.003$ & $0.002$ & $12.75$ & $5.060$ & $0.720$ &
$0.087$\\\hline
\textbf{DLo} & $0.150$ & $0.098$ & $1.830$ & $0.950$ & $-$ & $-$\\\hline
\textbf{Geo} & $0.582$ & $0.030$ & $-$ & $-$ & $-$ & $-$\\\hline
\textbf{DLi} & $0.436$ & $0.026$ & $-$ & $-$ & $-$ & $-$\\\hline
\textbf{P} & $1.390$ & $0.112$ & $-$ & $-$ & $-$ & $-$\\\hline
\textbf{DR} & $0.900$ & $0.009$ & $-$ & $-$ & $-$ & $-$\\\hline\hline
\end{tabular}
\end{align*}

\begin{align*}
& \text{Table 8. The goodness of fit tests for data set II. }\\
&
\begin{tabular}
[c]{|c|c|c|c|c|c|c|c|c|c|}\hline\hline
$\mathbf{X}$ & \textbf{OF} & \multicolumn{8}{|c|}{\textbf{EF}}\\\cline{3-10}
&  & \textbf{EDLi} & \textbf{DW} & \textbf{DB-XII} & \textbf{DLo} &
\textbf{Geo} & \textbf{DLi} & \textbf{P} & \textbf{DR}\\\hline\hline
\textbf{0} & $65$ & $64.97$ & $63.64$ & $63.32$ & $61.89$ & $45.98$ & $40.25$
& $27.42$ & $11$\\
\textbf{1} & $14$ & $14.39$ & $17.45$ & $18.19$ & $21.01$ & $26.76$ & $29.83$
& $38.08$ & $26.83$\\
\textbf{2} & $10$ & $9.01$ & $9.3$ & $9.29$ & $9.65$ & $15.57$ & $18.36$ &
$26.47$ & $29.55$\\
\textbf{3} & $6$ & $6.14$ & $5.68$ & $5.49$ & $5.24$ & $9.06$ & $10.35$ &
$12.26$ & $22.23$\\
\textbf{4} & $4$ & $4.33$ & $3.73$ & $3.52$ & $3.17$ & $5.28$ & $5.53$ &
$4.26$ & $12.49$\\
\textbf{5} & $2$ & $3.10$ & $2.56$ & $2.39$ & $2.06$ & $3.07$ & $2.86$ &
$1.18$ & $5.42$\\
\textbf{6} & $2$ & $2.24$ & $1.82$ & $1.69$ & $1.42$ & $1.79$ & $1.44$ &
$0.27$ & $1.85$\\
\textbf{7} & $2$ & $1.62$ & $1.32$ & $1.23$ & $1.02$ & $1.04$ & $0.71$ &
$0.05$ & $0.52$\\
\textbf{8} & $1$ & $1.18$ & $0.98$ & $0.92$ & $0.76$ & $0.61$ & $0.35$ &
$0.01$ & $0.11$\\
\textbf{9} & $1$ & $0.85$ & $0.74$ & $0.70$ & $0.58$ & $0.35$ & $0.17$ & $0$ &
$0.02$\\
\textbf{10} & $1$ & $0.62$ & $0.57$ & $0.55$ & $0.46$ & $0.21$ & $0.08$ & $0$
& $0$\\
\textbf{11} & $2$ & $1.55$ & $2.21$ & $2.71$ & $2.74$ & $0.28$ & $0.07$ & $0$
& $0$\\\hline
\textbf{Total} & $110$ & $110$ & $110$ & $110$ & $110$ & $110$ & $110$ & $110
$ & $110$\\\hline
$\mathbf{-L}$ &  & $166.9$ & $167.9$ & $168.8$ & $170.5$ & $178.8$ & $189.1$ &
$246.2$ & $277.8$\\
\textbf{AIC} &  & $337.9$ & $339.9$ & $343.5$ & $344.9$ & $359.5$ & $380.2$ &
$494.4$ & $557.6$\\
\textbf{CAIC} &  & $338.0$ & $340.1$ & $343.8$ & $345.1$ & $359.6$ & $380.3$ &
$494.5$ & $557.6$\\
\textbf{BIC} &  & $343.3$ & $345.4$ & $351.6$ & $350.4$ & $362.2$ & $382.9$ &
$497.1$ & $560.3$\\
\textbf{HQIC} &  & $340.1$ & $342.2$ & $346.8$ & $347.2$ & $360.6$ & $381.3$ &
$495.5$ & $558.7$\\\hline
$\mathbf{\chi}^{2}$ &  & $0.507$ & $1.04$ & $2.469$ & $3.316$ & $22.84$ &
$43.48$ & $294.1$ & $321.1$\\
\textbf{d.f} &  & $3$ & $3$ & $3$ & $3$ & $4$ & $4$ & $4$ & $4$\\
\textbf{P.value} &  & $0.917$ & $0.792$ & $0.480$ & $0.345$ & $<0.01$ &
$<0.01$ & $<0.01$ & $<0.01$\\\hline\hline
\end{tabular}
\end{align*}
From Table 8, it is clear that the EDLiD is the best distribution among all
tested models. Figure 6 shows the fitted PMFs for data set II, which support
the results in Table 8.%

\[%
{\parbox[b]{5.047in}{\begin{center}
\includegraphics[
height=4.0421in,
width=5.047in
]%
{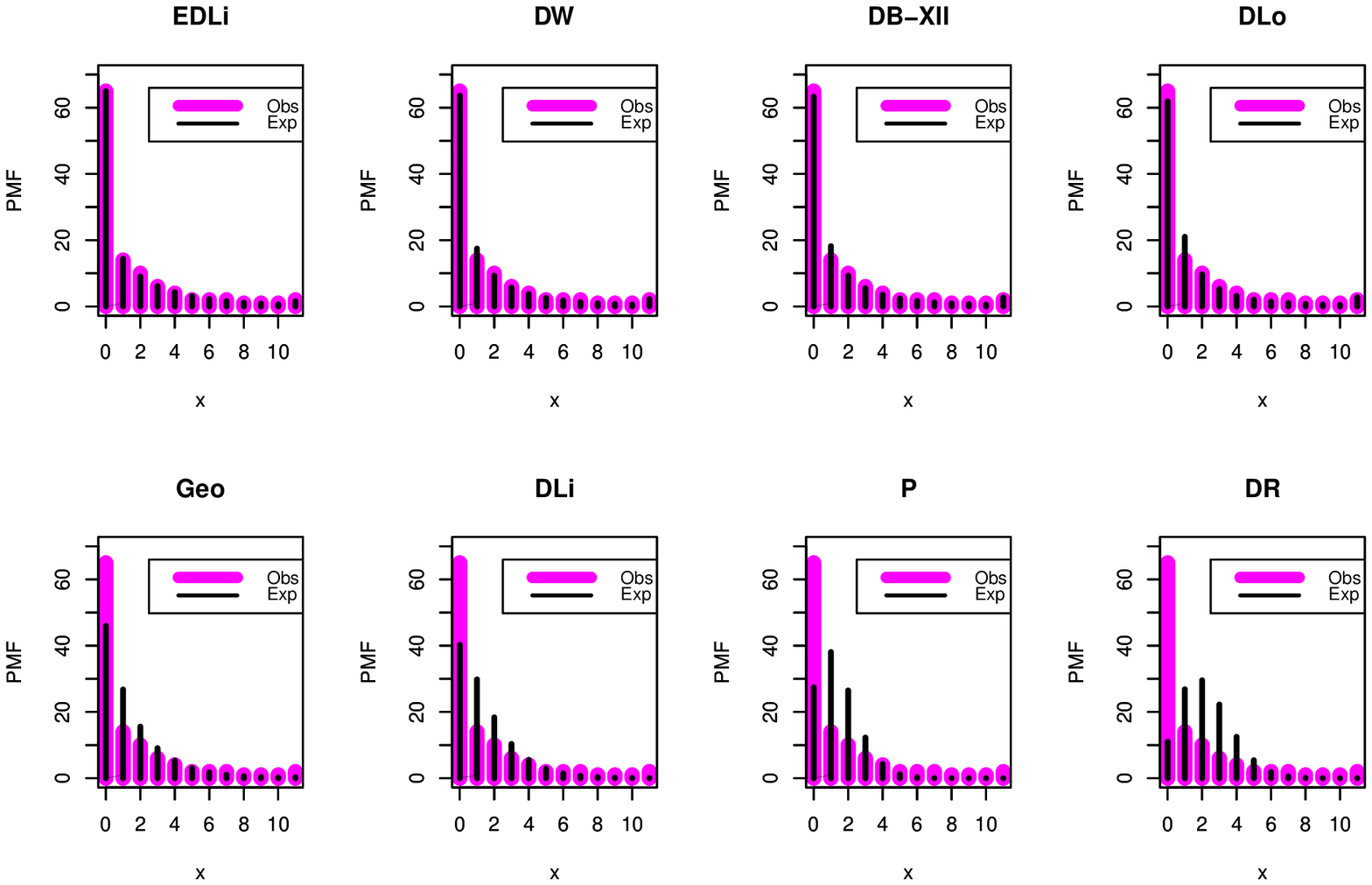}%
\\
Figure 6. The fitted PMFs for data set II.
\end{center}}}%
\]

\section{Conclusions}

A two-parameter EDLiD has been proposed. Its various distributional properties
have been discussed. It was found that the proposed distribution has a simple
structure, is more flexible and has a longer tail than the DLiD and other
discrete distributions in modeling data from different fields. In the future,
we will discuss the bivariate and multivariate extensions of this distribution.


\begin{thebibliography}{99}                                                                                               %
\bibitem {}Alamatsaz, M., Dey, H., Dey, S., Harandi, T., and Shams, S.,
(2016). Discrete generalized Rayleigh distribution. Pakistan journal of
statistics, 32(1), 1-20.

\bibitem {}Altun, G., Alizadeh, M., Altun, E., and \"{O}zel, G., (2017). Odd
Burr Lindley distribution with properties and applications. Hacettepe journal
of statistics and mathematics, 46 (2), 255-276.

\bibitem {}Bakouch, H. S., Aghababaei, M., and Nadarajah, S., (2014). A new
discrete distribution. Statistics, 48(1), 200-240.

\bibitem {}Bakouch, H. S., Al-Zahrani, B. M., Al-Shomrani, A. A., Marchi, V.
A., and Louzada, F., (2012). An extended Lindley distribution. Journal of the
Korean statistical society, 41, 75-85.

\bibitem {}Bebbington, M., Lai, C. D., Wellington, M., and Zitikis, R.,
(2012). The discrete additive Weibull distribution: a bathtub-shaped hazard
for discontinuous failure data. Reliability engineering and system safety,
106, 37-44.

\bibitem {}Calder\'{\i}n-Ojeda, E., and G\'{o}mez-D\'{e}niz, E., (2013). An
extension of the discrete Lindley distribution with applications. Journal of
the Korean statistical society, 42, 371-373.

\bibitem {}Chan, S., Riley, P. R., Price, K. L., McElduff, F., and Winyard, P.
J., (2009). Corticosteroid-induced kidney dysmorphogenesis is associated with
deregulated expression of known cystogenic molecules, as well as indian
hedgehog. American journal of physiology-renal physiology, 298(2), 346-356.

\bibitem {}Chandrakant, K., Yogesh, M. T., and Manoj, K. R., (2017). On a
discrete analogue of linear failure rate distribution. American journal of
mathematical and management sciences, 36(3), 229-246.

\bibitem {}Consul, P.C., and Jain, G. G., (1973). A generalization of the
Poisson distribution. Technometrics, 15(4), 791-799.

\bibitem {}Elbatal, I., Diab, L. S., and Elgarhy, M., (2016). Exponentiated
quasi Lindley distribution. International journal of reliability and
applications, 17(1), 1-19.

\bibitem {}Ghitany, M. E., Al-Mutairi, D. K., and Nadarajah, S., (2008a).
Zero-truncated Poisson Lindley distribution and its application. Mathematics
and computers in simulation, 79(3), 279-287.

\bibitem {}Ghitany, M. E., Al-Mutairi, D. K., Balakrishhnan, N., and Al-Enezi,
L. J., (2013). Power Lindley distribution and associated inference.
Computational statistics and data analysis, 64, 20-33.

\bibitem {}Ghitany, M. E., Alqallaf, F., Al-Mutairi, D. K., and Husain, H. A.,
(2011). A two-parameter weighted Lindley distribution and its applications to
survival data. Mathematics and computers in simulation, 81(6), 1190-1201.

\bibitem {}Ghitany, M. E., and Al-Mutairi, D. K., (2009). Estimation methods
for the discrete Poisson Lindley distribution. Journal of statistical
computation and simulation, 79(1), 1-9.

\bibitem {}Ghitany, M. E., Atieh, B., and Nadarajah, S., (2008b). Lindley
distribution and its application. Mathematics and computers in simulation,
78(4), 493-506.

\bibitem {}Goliforushani, S., and Asadi, M., (2008). On the discrete mean past
lifetime. Metrika, 68, 209-217.

\bibitem {}G\'{o}mez-D\'{e}niz, (2010). Another generalization of the
geometric distribution. Test, 19(2), 399-415.

\bibitem {}G\'{o}mez-D\'{e}niz, E., and Calder\'{\i}n-Ojeda, E., (2011). The
discrete Lindley distribution: properties and applications. Journal of
statistical computation and simulation, 81(11), 1405--1416.

\bibitem {}Hosking, J. R., (1990). L-moments: analysis and estimation of
distributions using linear combinations of order statistics. Journal of the
Royal statistical society, 52(B), 105-124.

\bibitem {}Inusah, S., and Kozubowski, T. J., (2006). A discrete analogue of
the Laplce distribution. Journal of statistical planning and infernce, 136, 1090-1102.

\bibitem {}Jehhan, A., Mohamed, I. , Eliwa, M. S., Al-mualim, S., and Yousof,
H. M., (2018). The two-parameter odd Lindley Weibull lifetime model with
properties and applications. International journal of statistics and
probability, 7, (4), 57-68.

\bibitem {}Jodr\'{a}, P., (2010). Computer generation of random variables with
Lindley or Poisson-Lindley distribution via the Lambert W function.
Mathematics and computers in simulation, 81(4), 851-859.

\bibitem {}Krishna, H., and Pundir, P. S., (2009). Discrete Burr and discrete
Pareto distributions. Statistical methodology, 6, 177-188.

\bibitem {}Kus, C., Akdogan,Y., Asgharzadeh, A., Kinaci, I., and Karakaya, K.,
(2018). Binomial-discrete Lindley distribution. Communications faculty of
sciences university of Ankara series A1-mathematics and statistics, 68(1), 401-411.

\bibitem {}Lehmann, E. L., (1952). The power of rank tests. Annals of
mathematical statistics, 24, 23-43.

\bibitem {}Lindley, D. V., (1958). Fiducial distributions and Bayes theorem.
Journal of the Royal statistical society, 20(B), 102-107.

\bibitem {}Liyanage, W., and Pararai, M., (2014). A generalized power Lindley
distribution with applications. Asian journal of mathematics and applications,
18, 1-23.

\bibitem {}Mahmoud, E., (2018). Logarithmic inverse Lindley distribution:
model, properties and applications. Journal of King Saud university science.
To appear.

\bibitem {}Mahmoudi, E., and Zakerzadeh, H., (2010). Generalized
Poisson-Lindley distribution. Communications in statistics: theory and
methods, 39(10), 1785-1798.

\bibitem {}Merovci, F., (2013). Transmuted Lindley distribution. International
journal of open problems in computer science and mathematics, 6(2), 63-72.

\bibitem {}Merovci, F., and Elbatal, I., (2014). Transmuted Lindley-geometric
distribution and its applications. Journal of statistics applications and
probability, 3(1), 77-91.

\bibitem {}Merovci, F., and Sharma, V. K., (2014). The Beta-Lindley
distribution: properties and applications. Journal of applied mathematics, http://dx.doi.org/10.1155/2014/198951.

\bibitem {}Munindra, B., Krishna, R. S., and Junali, H., (2015). A study on
two parameter discrete quasi Lindley distribution and its derived
distributions. International journal of mathematical archive, 6(12), 149-156.

\bibitem {}Nadarajah, S., Bakouch, H. S., and Tahmasbi, R., (2011). A
generalized Lindley distribution. Sankhya B, 73, 331-359.

\bibitem {}Nekoukhou, V., Alamatsaz, M. H., and Bidram, H., (2013). Discrete
generalized exponential distribution of a second type. Statistics, 47 (4), 876-887.

\bibitem {}Nedjar, S., and Zeghdoudi, H., (2016). On gamma Lindley
distribution: proprieties and simulations. Journal of computational and
applied mathematics, 298, 167-174.

\bibitem {}\"{O}zel, G., Alizadeh, M., Cakmakyapan, S., Hamedani, G., Ortega,
E. M., and Cancho, G., (2017). The odd log-logistic Lindley Poisson model for
lifetime data. Communications in statistics: simulation and computation,
46(8), 6513-6537.

\bibitem {}Pararai, M., Warahena-Liyanage, G., and Oluyede, B. O., (2015). A
new class of generalized power Lindley distribution with applications to
lifetime data. Theoretical mathematics and applications, 5, 53-96.

\bibitem {}Roy, D., (2003). The discrete normal distribution. Communications
in statistics: Theory and methods, 32, 1871-1883.

\bibitem {}Roy, D., (2004). Discrete Rayleigh distribution. IEEE transactions
on reliability, 53(2), 255-260.

\bibitem {}Sankaran, M., (1970). The discrete Poisson-Lindley distribution.
Biometrics, 26(1), 145-149.

\bibitem {}Shanker, R., and Mishra, A., (2013). A quasi Lindley distribution.
African journal of mathematics, 6(4), 64-71.

\bibitem {}Shanker, R., Sharma, S., and Shanker, R., (2013). A two-parameter
Lindley distribution for modeling waiting and survival times data. Applied
mathematics, 4, 363-368.

\bibitem {}Sharma, V., Singh, S., Singh, U., and Agiwal, V., (2015). The
inverse Lindley distribution: a stress-strength reliability model with
applications to head and neck cancer data. Journal of industrial and
production engineering, 32(3), 162-173.

\bibitem {}Tanka, R. A., and Srivastava, R. S., (2014). Size-biased discrete
two parameter Poisson-Lindley distribution for modeling and waiting survival
times data. Journal of mathematics, 10(1), 39-45.

\bibitem {}Vahid, N., and Hamid, B., (2015). The exponentiated discrete
Weibull distribution. Sort, 39 (1), 127-146.

\bibitem {}Zakerzadeh, H., and Dolati, A., (2009). Generalized lindley
distribution. Journal of mathematical extension, 3(2), 13-25.

\bibitem {}Zeghdoudi, H., and Nedjar, S., (2016). Gamma Lindley distribution
and its application. Journal of applied probability and statistics, 11(1), 129-138.

\bibitem {}Zeghdoudi, H., and Nedjar, S., (2017). On Poisson pseudo Lindley
distribution: Properties and applications. Journal of probability and
statistical science, 15(1), 19-28.
\end{thebibliography}
\end{document}